\documentclass[12pt,a4paper,twoside]{article}
\usepackage{amsmath,amsfonts,amssymb,a4,enumerate,epsfig,array,theorem,psfrag}

\newtheorem{theo}{Theorem}

\newtheorem{prop}{Proposition}[section]

\newcommand{\Xum}{{H_p^1}}
\newcommand{\Xdois}{{H^2_D}}

\newcommand{\ZZ}{{\mathbb{Z}}}
\newcommand{\RR}{{\mathbb{R}}}
\newcommand{\CC}{{\mathbb{C}}}
\newcommand{\Ss}{{\mathbb{S}}}
\newcommand{\BB}{{\mathbb{B}}}

\newcommand{\HH}{{\mathbb{H}}}
\newcommand{\hh}{{\mathbf{h}}}

\newcommand{\menos}{\setminus}
\newcommand{\pV}{\partial V}

\begin{document}
\title{The topology of critical sets of \\ some ordinary differential operators}
\author{Nicolau C. Saldanha and Carlos Tomei}
\maketitle

\textit{\hfill Dedicated to Djairo Figueiredo, with affection and admiration.}

\begin{abstract}
We survey recent work of Burghelea, Malta and both authors on the
topology of critical sets of nonlinear ordinary differential operators.
For a generic nonlinearity $f$, the critical set of the
first order nonlinear operator $F_1(u)(t) = u'(t) + f(u(t))$
acting on the Sobolev space $H^1_p$ of periodic functions
is either empty or ambient diffeomorphic to a hyperplane.
For the second order operator $F_2(u)(t) = -u''(t) + f(u(t))$
on $H^2_D$ (Dirichlet boundary conditions),
the critical set is ambient diffeomorphic to
a union of isolated parallel hyperplanes.
For second order operators on $H^2_p$,
the critical set is not a Hilbert manifold
but is still contractible and admits a normal form.
The third order case is topologically far more complicated.
\end{abstract}

\medbreak

\smallbreak

{\noindent\bf MSC-class:}
Primary 34L30, 58B05, Secondary 34B15, 46T05.

\section{Introduction}
\label{section:intro}

We survey recent work of Burghelea, Malta and both authors on the
topology of critical sets $C$
of (nonlinear, ordinary) differential operators $F$.
Our approach is \textit{geometric}, in the sense that we study
the geometry and topology of $C$ and its image $F(C)$
with the purpose of understanding the function $F$.
Pioneering examples of this point of view for nonlinear differential
equations are the results on the Laplacian coupled to a special
nonlinearity by Ambrosetti and Prodi (\cite{AP}),
interpreted as properties of a global fold by Berger and Podolak (\cite{BP}).

We begin this paper with a two dimensional example in which
most of our claims can be followed visually.
We then consider differential operators of increasing difficulty.
For generic nonlinearities $f: \RR \to \RR$,
set $F_1(u) = u'+ f(u)$ on the Sobolev space
$H^1_p$ of periodic functions (\cite{MST1}), with critical set $C_1$.
We show that if $C_1$ is not empty then there is
a global diffeomorphism in the domain of $F_1$
converting $C_1$ into a hyperplane.
There are two parts in the argument.
First, we prove that the homotopy groups of $C_1$ are trivial.
We then use the fact that, under very general conditions,
(weak) homotopically equivalent
infinite dimensional Hilbert manifolds are actually diffeomorphic.
For some classes of nonlinearities, we achieve a global
normal form for $F_1$. 

We proceed to the second order operator 
$F_{2,D}(u) = -u'' + f(u)$ acting on $H^2_D([0,\pi])$,
the Sobolev space of functions satisfying Dirichlet boundary
conditions (\cite{BST}).
The critical set $C_{2,D}$ now splits into connected components
$C_{2,D,m}$, one for each positive integer value of $m$
for which $-m^2$ is in the interior of the image of $f'$.
Again, there exists a diffeomorphism taking each $C_{2,D,m}$
to a hyperplane.
The proof of the triviality of homotopy groups of $C_{2,D,m}$
requires a very different strategy than the first order case.

Second order operators on $H^2_p$ (the space of periodic
functions) are more complicated, in the sense that the critical
set is not a manifold, but with suitable hypothesis
is still contractible and admits a normal form (\cite{BST2}).
In brief, the critical set looks like a union of a hyperplane
and infinitely many cones, with nonregular points at the vertices.
Nonregular points in the critical set correspond to
potentials $h(t) = f'(u(t))$, $u \in H^2_p$, for which the kernel
of the linear operator $v \mapsto -v'' + h(t) v$ is a subspace of
dimension $2$ of $H^2_p$. The diffeomorphism of $H^2_p$
through which the critical set achieves its normal form
takes the set of nonregular critical points to a disjoint union
of linear subspaces of codimension $3$.

The third order case is yet more complicated:
consider the set $C^\ast_{3,p}$ of pairs $(h_0,h_1)$
for which the kernel
of the linear operator $u \mapsto u''' + h_1 u'+ h_0 u$ is
a subspace of $H^3_p$ of dimension $3$.
In the final section we show that
the space $C^\ast_3$ is homotopically equivalent to $X_I$,
the set of closed locally convex curves in the sphere $\Ss^2$
with a prescribed initial base point and direction.
It turns out that the space $X_I$ has
three connected components (\cite{Little}),
two of them having a rich algebraic topological structure,
not equivalent to any finitely generated CW-complex
(\cite{S}, \cite{SK}, \cite{Shapiro2}, \cite{ShapiroM}).

The authors acknowledge the support of CNPq, Faperj and Finep.

\vfil\eject
\section{A finite dimensional example}

Consider
\begin{align} F: \CC &\to \CC \notag\\
z &\mapsto z^7 + 5\bar{z}^4 + z \notag \end{align} which is
clearly a smooth (but not analytic) function from the plane
$\RR^2$ to itself.
In this section, we identify $\CC$ and $\RR^2$ indiscriminately.
How many solutions has the equation $F(z)=0$?

Our approach to the question (\cite{MST0}, \cite{ST})
breaks into a few steps. First, we
compute the critical set $C$ of $F$ by
searching for points in the plane
in which the Jacobian $DF$ is not invertible.
Some numerical analysis reveals that $C$ consists of two simple curves
$C_i$ and $C_o$ bounding disks $(0,0) \in D_i \subset D_o$.
A finer inspection verifies that, as expected from Whitney's
classical theorems on planar singularity theory (\cite{W}), the
generic critical point is a \textit{fold},
i.e., after local changes of variables
the function takes the form $(x,y) \mapsto (x,y^2)$ close to the origin.
Also, nonfold points are \textit{cusps}, 
with local normal form $(x,y) \mapsto (x, y^3 + xy)$.
Again, with the help of some computation, one finds out
that the curves $C_i$ and $C_o$ have five and eleven cusps,
respectively. What may be more informative is the geometry of the
images $F(C_i)$ and $F(C_o)$, presented in the picture below.

\begin{figure}[ht]
\begin{center}
\epsfig{height=60mm,file=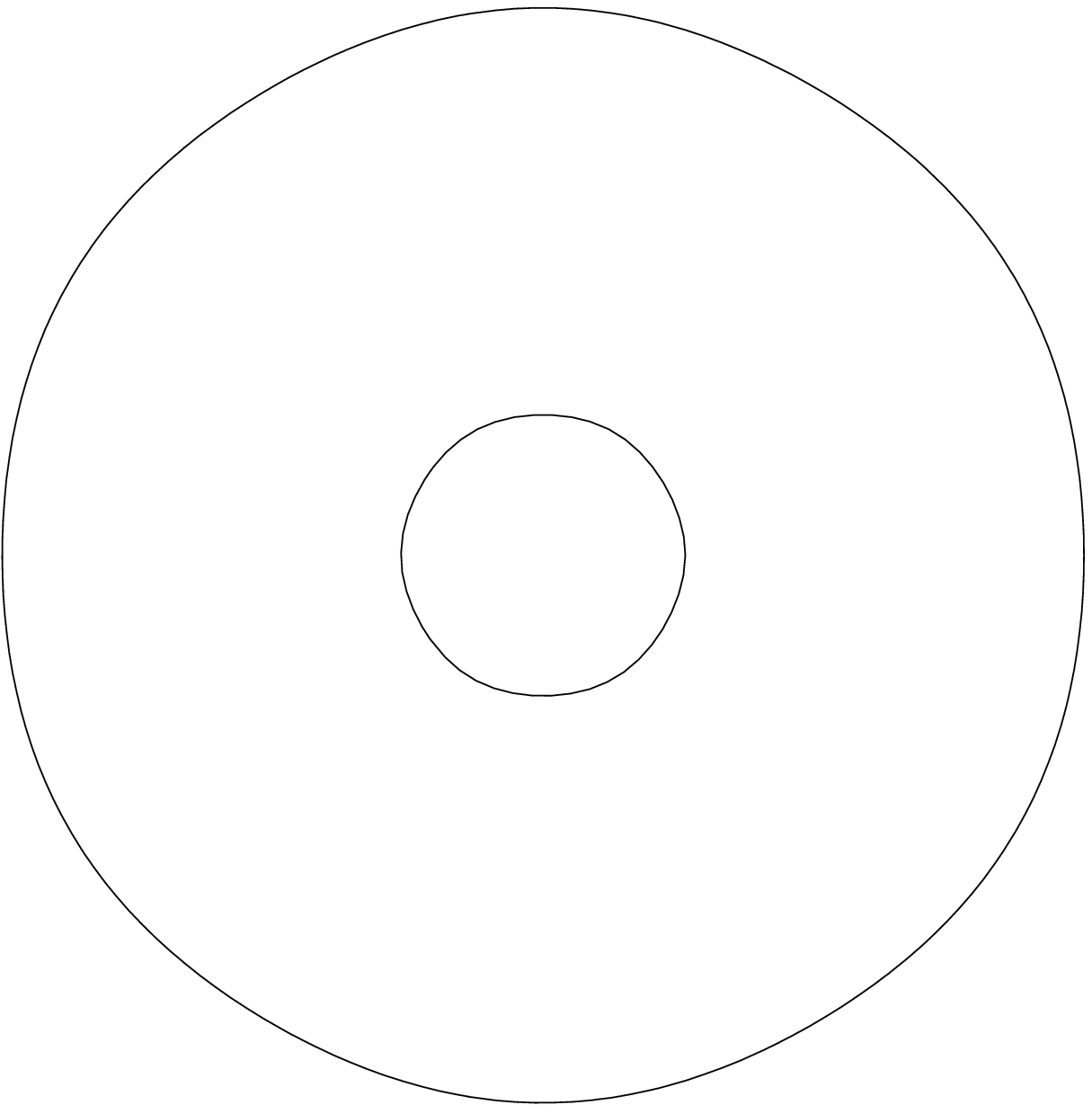} \quad
\epsfig{height=60mm,file=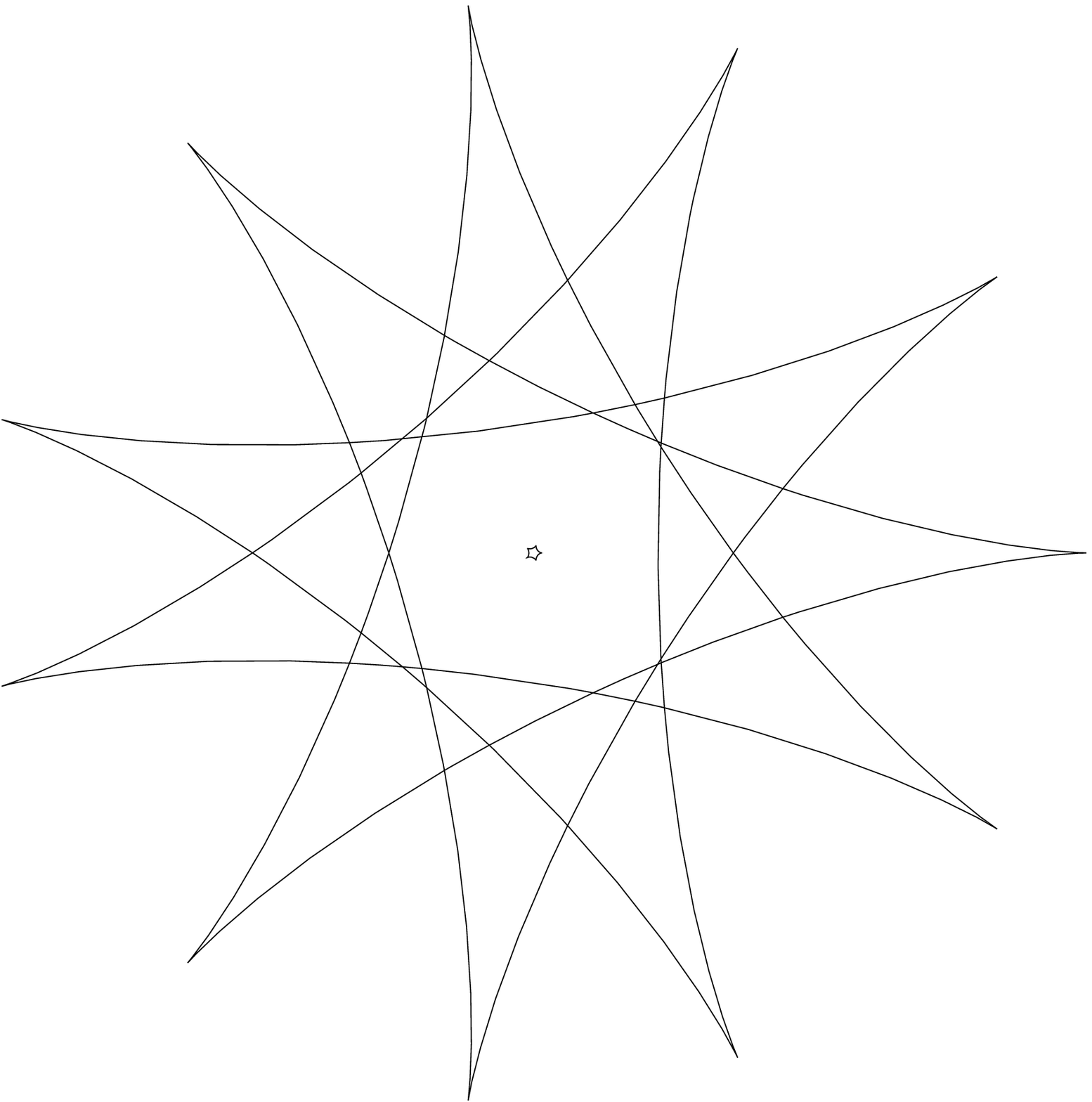}
\end{center}
\caption{The critical set of the function $F(z) = z^7 + 5\bar{z}^4 + z$
and its image.}
\label{fig:z74}
\end{figure}

We now count preimages. Simple estimates suffice to prove that
$F$ is proper, and that, for sufficiently large
$(w_1,w_2) \in \RR^2$, the equation $F(x,y) = (w_1,w_2)$ has seven solutions.
From properness, points in the same connected component of the
complement of $F(C)$ have the same number of pre-images.
Now, consider two points $p$ and $q$ in components sharing a
boundary arc of (images of) folds. From the normal form of a
function at a fold point, the number of preimages of $p$ and $q$
under $F$ differ by two.
The {\it sense of folding} indicates that
the number of preimages {\it increases} by two
whenever one gets closer to the origin when trespassing
an arc of images of folds.
Adding up, from the knowledge that points in the unbounded
connected component of the complement of $F(C)$ have seven
preimages, we learn the number of preimages of each component.
Thus, the origin has 17 preimages. In a nutshell, the number
of preimages of a point under a function $F$ may become large when
the image of the critical set $C$ intersects itself abundantly.


In the following sections, we shall consider functions between
separable infinite dimensional Hilbert spaces.
In many cases, the critical set turns out to be surprisingly simple.

\section{The first order operator}
\label{section:topo}

We consider the differential equation
\[ u'(t) + f(u(t)) = g(t), \qquad u(0) = u(\pi) \]
where the unknown $u$ is a real function
and $f: \RR \to \RR$ is a smooth nonlinearity.
In the spirit of the example in the previous section, we define the operator
\begin{align} F_1: \Xum &\to L^2 \notag\\
u &\mapsto u'+ f(u)  \notag
\end{align}
where $\Xum = H_p^1([0,1];\RR) = H^1(\Ss^1;\RR)$
is the Sobolev space of periodic functions
with weak derivative in $L^2 = L^2([0,1];\RR)$.
It is easy to verify that the differential
\[ DF_1(u) v = v'+ f'(u)v \]
is a Fredholm operator of index $0$
and therefore $u$ belongs to the critical set $C_1 \subset \Xum$
of $F_1$ if and only if the equation
\[ v'(t) + f'(u(t)) v(t) = 0, \qquad v(0) = v(1) \]
admits a nontrivial solution $v$, i.e.,
\[ C_1 = \{u \in \Xum  \;|\;  \phi_1(u) = 0 \}, \qquad
\phi_1(u) = \int_0^1 f'(u(t)) dt. \]
An equivalent spectral interpretation for $C_1$ is
\[ C_1 = \left\{ u \in \Xum \;|\;
0 \in \sigma(DF_1(u)) \right\}. \]
It is not hard to see (\cite{MST1}) that $DF_1(u)$ has a unique
real eigenvalue.
We assume that $f''$ has isolated roots which are
distinct from the roots of $f'$:
this implies that $0$ is a regular value of $\phi_1$
and therefore that $C_1$ is a smooth hypersurface in $\Xum$.
An $\HH$-manifold is a manifold modeled on the separable infinite
dimensional Hilbert space $\HH$: $C_1$ is an $\HH$-manifold.
With these hypothesis, the topology of $C_1$ is trivial (\cite{MST1}):

\begin{theo}\label{theo:A}
Assume $C_1$ to be nonempty.
Then $C_1$ is path connected and contractible.
Furthermore, there is a diffeomorphism
from $\Xum$ to itself taking $C_1$ to a hyperplane.
\end{theo}

A natural finite dimensional counterpart of this theorem is false.
Indeed, let
\[ C_1^n = \{u \in \RR^n  \;|\;  \phi_1^n(u) = 0 \}, \qquad
\phi_1^n(u) = \sum_k f'(u_k). \]
For $f(x) = x^3 - x$, $C_1^n$ is a sphere.
Thus, our theorem goes hand in hand with the well known facts
that in an infinite dimensional Hilbert space $\HH$,
the unit sphere is contractible
and there is a diffeomorphism of $\HH$ to itself
taking the unit sphere to a hyperplane.

\begin{figure}[ht]
\begin{center}
\psfrag{t}{$t$}
\psfrag{s}{$s$}
\epsfig{height=40mm,file=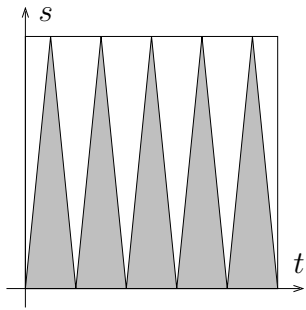}
\end{center}
\caption{The set $A$ (shaded).}
\label{fig:minhoca}
\end{figure}

The rest of this section outlines the proof of theorem \ref{theo:A}.
First, one proves (\cite{MST2}) that the homotopy groups
$\pi_k(C_1)$ are trivial.
We consider only $k = 0$, i.e., path connectedness of $C_1$,
the other cases being similar.
This is done first in the space $C^0$ (with the sup norm): more precisely, set
\[ C^0_1 = \left\{ u \in C^0([0,1]) | \phi_1(u) = 0 \right\}; \]
it is easier to build a homotopy in the sup norm, since it is weaker.
Take two functions $u_0(t)$ and $u_1(t)$ in $C^0_1$ so that
\[ \int_0^1 f(u_0(t)) dt = \int_0^1 f(u_1(t)) dt =0.\]
Let $s \in [0,1]$ be the parameter for the desired homotopy
$\hh(s,t)$: we want $\hh(0,t)=u_0(t), \hh(1,t)=u_1(t), \hh(s, \cdot) \in C^0_1$.
Define a discontinuous function $\hh_0: [0,1] \times [0,1] \to \RR$
by $\hh_0(s,t) = u_0(t)$ if $(s,t) \in A \subset [0,1] \times [0,1]$
and $\hh_0(s,t) = u_1(t)$ otherwise,
where $A$ looks like the set in figure \ref{fig:minhoca}:
for each $s$, $A_s = \{t \;|\; (s,t) \in A \}$ is a rather uniformly
distributed subset of $[0,1]$ of measure $1-s$.
For each $s$, $\phi_1(\hh_0(s,\cdot)) \approx 0$.
The function $H$ equals $\hh_0$ except on a set of very small measure,
where it is defined so as to obtain a continuous $\hh$
for which $\phi_1(\hh(s,\cdot)) = 0$.

In order to prove that the homotopy groups $\pi_k(C_1)$ are also trivial,
we use the following theorem (\cite{BST}) with $Y = \Xum$, $X = C^0([0,1])$,
$M = C_1^0$, $N = C_1$ and $i$ the inclusion.

\begin{theo}\label{theo:C}
Let $X$ and $Y$ be separable Banach spaces. Supose
$i: Y \to X$ is a bounded, injective map with dense image and $M \subset X$
a smooth, closed submanifold of finite codimension.
Then $N = i^{-1}(M)$ is a smooth closed submanifold of $Y$, and the
restrictions $i: Y-N \to X-M$ and $i: (Y,N) \to (X,M)$ are
homotopy equivalences.
\end{theo}

Since $C_1$, being an $\HH$-manifold, is homotopy equivalent to a CW-complex,
weak contractibility implies contractibility and $C_1$ is contractible
(and, by theorem \ref{theo:C}, $C_1^0$ is contractible too).
The connected components of the complement of $C_1$ are also contractible.
Indeed, let $h: \Ss^k \to \Xum - C_1$ be a continuous map;
since $\Xum$ is obviously contractible, $h$ can be extended
to $\hh: \BB^{k+1} \to \Xum$.
Since $C_1$ is contractible, $\hh$ can be redefined to take
each connected component of $\hh^{-1}(\Xum - C_1)$
not touching $\Ss^k$ to $C_1$. A normal vector to $C_1$
can then be used to push the image of $\hh$ away from $C_1$,
yielding $\hh: \BB^{k+1} \to \Xum - C_1$.
Thus, the homotopy groups $\pi_k(\Xum - C_1)$, $k \ge 1$, are trivial
and therefore the connected components of $\Xum - C_1$ are contractible.
To complete the proof of theorem \ref{theo:A},
we use the following result (\cite{BST}):

\begin{theo}\label{theo:B}
Suppose $f:(V_1, \pV_1) \to (V_2, \pV_2)$ is a smooth homotopy
equivalence of $\HH$-manifolds with boundary, $K_2 \subset V_2
\menos \pV_2$ a closed submanifold of finite codimension and $K_1
= f^{-1}(K_2)$. Suppose also that $f$ is transversal to $K_2$ and
the maps $f:K_1 \to K_2$ and $f: V_1 - K_1 \to V_2 - K_2$
are homotopy equivalences. Then there exists a diffeomorphism
between $(V_1; \pV_1, K_1)$ and $(V_2; \pV_2, K_2)$ which is homotopic to
$f$ as maps of triples.
\end{theo}

Set $V_1 = \Xum$ and
consider the functional $\phi_1: \Xum \to \RR$ defined above.
Write $V_2 = L^2([0,1]) = \langle 1 \rangle \oplus \langle 1
\rangle^{\perp} $, where $\langle 1 \rangle \cong \RR$ is the
vector space of constant functions and $\langle 1
\rangle^{\perp}$, the space of functions $v$ of average $\bar{v}$
equal to zero, is a hyperplane in $V_2$. Set $f : V_1 \to V_2$
to be $f(u) = (\phi_1(u), 0)$. Set $K_2 = \langle 1 \rangle^\perp$
so that the critical set $C_1$ equals
$K_1 = f^{-1}(K_2)$, and the genericity condition
on $f$ ensures that $C_1 = K_1$ is a hypersurface of
$V_1$ and $f$ is transversal to $K_2$.
Notice that $\pV_1 = \pV_2 =\emptyset$.
Since $K_1 = C_1$ is contractible, we obtain that $f$ is a
homotopy equivalence between $K_1$ and $K_2$.
Similarly, since the connected components of $V_1 - K_1$
are contractible, $f: V_1 - K_1 \to V_2 - K_2$ is a homotopy equivalence
and Theorem \ref{theo:A} follows.

The proofs of theorems \ref{theo:C} and \ref{theo:B} are rather technical
and shall not be discussed here.
A simple spinoff is the following
density result which may be of independent interest.

\begin{theo}
\label{theo:dense}
Let $X$ be a Banach space, $V \subset X$ a dense (linear) subspace
and $M \subset X$ a finite codimension submanifold. Then the
intersection $V \cap M$ is dense in $M$.
\end{theo}

\section{More on $C_1$ and $F_1(C_1)$}

For the first order differential operator $F_1$, generically, it
is easy to prove that the critical $C_1$ stratifies by complexity
of the singularities, which are infinite dimensional counterparts
of the familiar Morin singularities (\cite{Morin}, \cite{MST1}).
In a nutshell, the Morin singularity is the generic situation for
a germ from $(\RR^n,0)$ to itself whose Jacobian at 0 has a one
dimensional kernel.  There is essentially one \textit{elementary}
Morin singularity in each dimension, with very simple normal
forms: the fold $x \mapsto x^2$, the cusp $(x,y) \mapsto (x, y^3
+xy)$, the swallowtail $(x,y,z) \mapsto (x, y,z^4 + y z^2+xz)$,
the butterfly $(x,y,z,t) \mapsto (x, y,z,t^5 + zt^3+yt^2+xt)$,\dots
One should take into account Morin singularities arising by
taking cartesian products: the fold $x \mapsto x^2$ embeds into
the (non-elementary) fold $(x,y) \mapsto (x^2,y)$. Even in
infinite dimensions, Morin singularites are always cartesian
products of an elementary singularity (for which one defines a
natural \textit{order}, given by the degree of the top polynomial
of its normal form) and the identity mapping with the appropriate
codimension (\cite{MST1}; also \cite{CDT} for infinite dimensional
cusps).

What is the topology of the finer strata? For the first order
differential operator $F_1$ with a generic nonlinearity $f$, the
subset $\Sigma_2$ of critical points which are nonfold points
again has trivial topology. The proof of this statement begins
with the following result:

\begin{theo}
\label{theo:Sigma2} For a generic smooth nonlinearity $f: \RR \to \RR$,
there is a diffeomorphism from $\Xum$ to itself taking the sets $C_1$ and
$\Sigma_2$ of $F_1$ respectively to the zero level of the
functional $\phi_1$ and to the zero level of the pair
$(\phi_1,\phi_{1,2})$, where
\[ \phi_{1,2}(u) = \int_0^1 f''(u(t)) dt. \]
\end{theo}

Consider now a vector valued extension of the argument used to
prove that the homotopy groups of $C_1$ are trivial (\cite{MST2}).
Let $M$ be a compact, finite dimensional manifold with a smooth
Riemannian metric inducing a normalized measure $\mu$, so that
$\mu(M)=1$. Let $V$ a separable Banach space of continuous real
valued functions on $M$, which is closed under multiplication by
functions of $C^\infty(M)$ (here, multiplication in $C^\infty(M)
\times V \to V$ is continuous). For a smooth function $f_n: M
\times \RR \to \RR^n$ define $N_n: V \to \RR^n$ as the average of
the Nemytski{\u\i} operator associated to $f$,
\[ N_n(v) = \int_M f_n(m,v(m)) d\mu.\]
Let $\Pi_k : \RR^n \to \RR^k$ be the projection on the first $k$
coordinates. We say that $0$ is a \textit{strong regular value} of
$N_n$ if $0$ is a regular value of each composition $N_k = \Pi_k
\circ F_n$. Finally, let $Z_k = N_k^{-1}(0)$: if $0$ is a strong
regular value of $N_n$, the sets $Z_k$ are nested submanifolds of
$V$ of codimension $k$.

\begin{theo}
\label{theo:Nemi}
Suppose $0$ is a strong regular value of $N_n$ as
above. Then the levels $Z_k$ are contractible. Moreover, there is
a global homeomorphism $\Psi: V \to V$ taking each $Z_k$ to a
closed linear subspace of $V$ of codimension $k$. If $V$ is a
Hilbert space, $\psi$ can be taken to be a diffeomorphism.
\end{theo}

Thus, for a generic nonlinearity, the zero levels of  $\phi_1$ and
$(\phi_1,\phi_{1,2})$ are ambient diffeomorphic to nested
subspaces of codimension 1 and 2, and we obtain the topological
triviality of $C_1$ (again) and of $\Sigma_2$. Theorem
\ref{theo:Sigma2} is necessary because the set of nonfold points
is \textit{not} described naturally in terms of zero levels of
vector Nemytski{\u\i} operators. This is the reason
for which we do not know how to address topological properties
of higher order strata.

The image $F(C_1)$ clearly depends on the nonlinearity. We list a
few facts (\cite{MST1}), some of which had been known (\cite{MS}).

\begin{enumerate}
\item{Suppose $f$ is strictly increasing, onto the interval $(a,b)$.
Then $F_1$ is a diffeomorphism onto the strip
\[ \{ w \in L^2 \; | \; a < \langle w, 1 \rangle < b \}. \]}
\item{If $f$ is strictly convex, $F_1$ is a global fold. If $f$ is generic,
and all critical points of $F_1$ are folds, then $F_1$ is a global
fold. There are nonconvex linearities for which $F_1$ is a global
fold.}
\item{If $f$ is a generic polynomial of even degree, positive leading
coefficient, taking  negative values, then $F_1$ has cusps.}
\item{Suppose $f: \RR \to \RR$ is proper, with derivative $f'$ assuming both signs
and third derivative with isolated roots assuming only one sign.
Then $F_1$ is a global cusp, in the sense that diffeomorphisms on
$\Xum$ and $L^2$ convert $F_1$ to the normal form
\begin{align} \tilde{F_1}: \RR^2 \times \HH &\to \RR^2 \times \HH \notag\\
(x,y;v) &\mapsto (x,y^3+xy;v)  \notag
\end{align}}
\end{enumerate}
The last example establishes a conjecture by Cafagna and Donati
(\cite{CD}) on the global topology of operators associated to the
nonlinearity $f(x) = ax + b x^2 + c x ^{2k+1}$, $a \ge 0$, $a^2 +
b^2 > 0$ and $c <0$.

From the examples, the number of preimages of an operator
associated  to a generic nonlinearity given by a polynomial of
degree less than or equal to three is bounded by the degree of the
polynomial. This is false for polynomials of degree four. The
counterexample  in \cite{MST1} was obtained  by adjusting
coefficients of the nonlinearity in order to obtain a butterfly in
the critical set of $F_1$. Singularity theory combined with a
degree theoretic argument imply the existence of a point with six
preimages, five of which are near the butterfly. Both nonlinearity
and special point were computed numerically.

On a related note, Pugh conjectured that the equation
\[ u'(t) = a_k(t) (u(t))^k + \cdots + a_1(t) u(t) + a_0(t), \qquad
u(1) = u(0) \]
would have at most $k$ solutions. The conjecture was essentially
verified by Smale for $k = 2$ and $k = 3$ if $a_3(t) > 0$ for all $t$.
The general conjecture was proved false by Lins Neto (\cite{Lins}).
Our example has the special feature that $a_i$ is constant in $t$ for $i > 0$.

\section{The second order case, Dirichlet conditions}

Let $\Xdois = H^2_D([0,\pi],\RR)$ be the Sobolev space
of functions with second weak
derivative in $L^2 = L^2([0,\pi];\RR)$, satisfying Dirichlet boundary
conditions and for a smooth nonlinearity $f$ now consider the
operator
\begin{align}F_{2,D}: \Xdois &\to L^2 \notag\\
u &\mapsto -u''+ f(u) \notag
\end{align}
with differential given by $DF_{2,D}(u)v = -v'' + f'(u)v$.
Again, by Fredholm theory, the critical set $C_{2,D}$ of
$F_{2,D}$ consists of functions $u \in \Xdois$ for which $DF_{2,D}$
has zero in the spectrum.

For a more explicit description of
$C_{2,D}$, define the fundamental solution $v_1$,
\[ -v_1'' + f'(u)v_1 = 0, \qquad v_1(0)= 0, v_1'(0) = 1, \]
and consider the (continuously defined) argument $\omega(t)$ of
the vector $(v'_1(t),v(t))$ for which $\omega(0)=0$. Notice the
implicit dependence of $\omega$ on $v_1$, which in turn depends on $u$.
Then
\[ C_{2,D} = \{ u \in \Xdois \; | \;  \omega(\pi)/\pi \in \ZZ \}. \]
Thus, there must be different connected components of $C_{2,D}$
for each value of $\omega(\pi) = m\pi$.
Said differently, the fact that the eigenvalue $0$ of $DF_{2,D}(u)$
may be the $m$-th eigenvalue in the spectrum (counting
from below) induces different critical components of $F_{2,D}$.
Define the functional $\phi_{2,D}(u)$ as the value of the argument
at $\pi$ of the fundamental solution $v_1$ associated to $u$ and
consider subsets $C_{2,D,m}$ which partition $C_{2,D}$, given by the
$m\pi$-levels of $\phi_{2,D}$.
For a generic set of nonlinearities $f$, the set
$C_{2,D,m}$ is nonempty if and only if the number $-m^2$ belongs to the
interior of the image of $f'$:
in this case, $C_{2,D,m}$ is path connected and contractible.


We want to prove that arbitrary functions from a sphere $\Ss^k$
to $C_{2,D,m}$ are homotopic to constants.
Theorem \ref{theo:Nemi} is not the appropriate tool anymore,
among other reasons because the functional $\phi_{2,D}$ depends 
on the derivative $u'$.
From the topological theorems of section 3,
it suffices to control the deformation in the sup norm.
Contractibility then implies a global change of variables in
$\Xdois$ which flattens the connected components $C_{2,D,m}$ into hyperplanes.
The stratification of the critical set in Morin
singularities still applies, but there are no results at this
point about the global topology of the finer strata.

\begin{figure}[ht]
\begin{center}
\epsfig{height=55mm,file=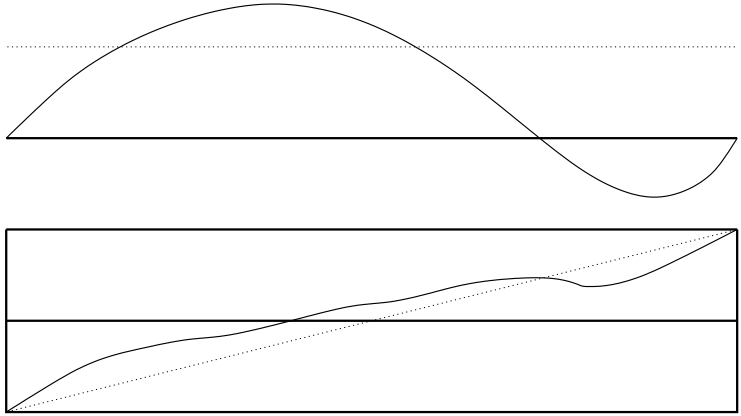}
\end{center}
\caption{Graphs of $u$ and $x_m$ and their $m$-arguments}
\label{fig:H0}
\end{figure}

The idea behind the construction of the homotopy
is to change $u$ so as to approach a function $u_\ast$
which is constant equal to a value $x_m$
with the property that $f'(x_m) = -m^2$
except for two small intervals at the endpoints.
Notice that the (fundamental) solution of
\[ -v'' -m^2 v = 0, \qquad v(0) = 0, v'(0) = 1 \]
is the function $v(t) = \sin(mt)/m$.
In general, define the \textit{$m$-argument} $\omega_m(t)$,
the argument of the vector $(v_1'(t),mv_1(t))$.
Since $\omega(t) = \omega_m(t)$ if $\omega(t) = k\pi$, $k \in \ZZ$,
the sets $C_{2,D,m}$ are $m\pi$-levels
for both $\omega(\pi)$ and $\omega_m(\pi)$.
Better still, the $m$-argument of $v(t) = \sin(mt)/m$
varies linearly from 0 to $m\pi$ for $t \in [0,\pi]$.
In figure \ref{fig:H0} we show an example of $u \in C_{2,D,m}$, $m = 2$,
and in dotted lines the constant value $x_m$;
below we show the $m$-argument $\omega_m$ of $u$
and the constant function $x_m$.

\begin{figure}[ht]
\begin{center}
\epsfig{height=55mm,file=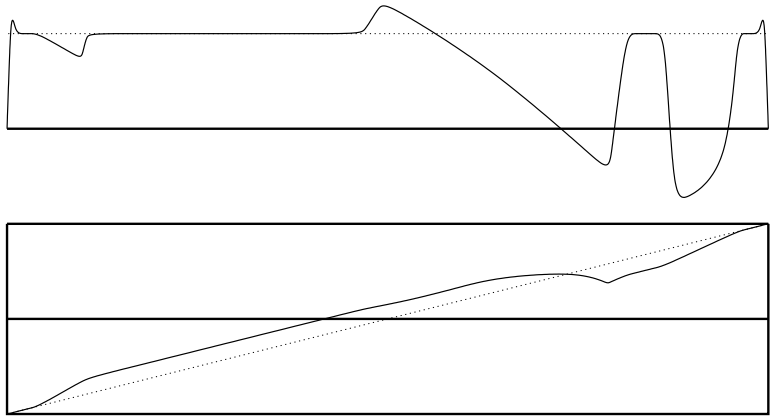}
\end{center}
\caption{The function $u$ gets squeezed}
\label{fig:H5}
\end{figure}

The homotopy squeezes the graph of $\omega_m$
between parallel walls advancing towards the line $y = mt$,
as shown in figure \ref{fig:H5}.
A corresponding $u$ is obtained by changing its original value
in the region of the domain over which the graph of $\omega_m$
has been squeezed---there, the new value of $u$ is $x_m$.
The value of $\omega_m(\pi)$ for this new $u$ equals $m \pi$,
but such $u$ is discontinuous and therefore not acceptable.
We make amends: the region where the graph of $\omega_m$
trespasses the wall by more than a prescribed tolerance is taken to $x_m$
and in the region where the graph of $\omega_m$ lies strictly
between the walls, $u$ is unchanged.
Hence, there is an open region in the domain
where $u$ assumes rather arbitrary values
in order to preserve its continuity.

Similar results were previously obtained for the special case when
$f$ is a convex nonlinearity (\cite{Ruf}, \cite{BT}). The trivial
topology of each connected component $C_{2,D,m}$ then follows from the
fact that each such set is a graph over the hyperplane in $\Xdois$
of functions orthogonal to $\sin t$.

\section{The second order periodic case}

Things get more complicated for the second order operator with periodic
boundary conditions. As is well known, the associated
linearizations do not necessarily have simple spectrum.
Critical sets now are \textit{not} submanifolds,
and one has to be careful about nonregular points.
Still, a description of the global geometry of the critical set is
possible (\cite{BST2}).

For any smooth nonlinearity $f: \RR \to \RR$
denote by $F_{2,p}$ the differential operator
\begin{align} F_{2,p}: H_p^2 &\to L^2 \notag\\
u &\mapsto -u'' + f(u) \notag \end{align}
where $H^2_p = H^2_p([0,2\pi]; \RR) = H^2(\Ss^1; \RR)$
and $L^2 = L^2([0,2\pi]; \RR) = L^2(\Ss^1; \RR)$.
We are interested in the critical set
$C_{2,p} \subset H^2_p$ of $F_{2,p}$.
Again, the differential of $F_{2,p}$
is the Fredholm linear operator $DF_{2,p}(u) v = -v'' + f'(u) v$ of index 0.

Let $\Sigma_0 \subset \RR^3$ be the plane $z = 0$ and,
for $n > 0$, let $\Sigma_n$ be the cone
\[ x^2 + y^2 = \tan^2 z, \quad
2\pi n - \frac{\pi}{2} < z < 2\pi n + \frac{\pi}{2} \]
and $\Sigma = \bigcup_{n \ge 0} \Sigma_n$.

\begin{theo}
\label{theo:periodic}
Let $f: \RR \to \RR$ be smooth function such that
$f'$ has isolated critical points and $f'$ is surjective.
Then the pair $(H^2_p, C_{2,p})$ is diffeomorphic to the pair
$(\RR^3 \times \HH, \Sigma \times \HH)$.
\end{theo}

The set $C^\ast_{2,p} \subset C_{2,p}$ of potentials $u$
for which {\it all} solutions $v$ of 
\[  - v''(t) + f'(u(t)) v(t) = 0 \]
are periodic is the set of nonregular points.
It follows from theorem \ref{theo:periodic} that $C^\ast_{2,p}$
is a submanifold of codimension $3$, taken by
the diffeomorphism mentioned in the theorem to
the set vertices of the cones, i.e.,
to points of the form $((0,0,2\pi n),\ast) \in \RR^3 \times \HH$.

The proof of theorem \ref{theo:periodic} requires the topological
study of the {\it monodromy map}. For $h \in L^2$,
let $v_1, v_2 \in H^2([0,2\pi]; \RR)$ be defined by
\[ v_i''(t) = h(t) v_i(t),\quad
v_1(0) = 1, \; v_1'(0) = 0, \; v_2(0) = 0, \; v_2'(0) = 1 \]
and define $\tilde\beta: [0,2\pi] \to G$ by
$\tilde\beta(0) = I$ and
\[ \tilde\beta(t) = \begin{pmatrix}
v_1(t) & v_1'(t) \\ v_2(t) & v_2'(t) \end{pmatrix} \]
where $G = \widetilde{SL(2,\RR)}$ is the universal cover of $SL(2,\RR)$.
Finally, define $\xi_\bullet: L^2 \to G$
by $\xi_\bullet(h) = \tilde\beta(2\pi)$.
The map $\xi_\bullet$ is not surjective:
its image is an open set $G_+ \subset G$ diffeomorphic to $\RR^3$.
The map $\xi_\bullet$ is in a sense a projection:

\begin{prop}
\label{prop:xibu}
There exists a smooth diffeomorphism $\Psi_k: G_+ \times \HH \to L^2$
such that $\xi_\bullet \circ \Psi_k$ is the projection
on the first coordinate.
\end{prop}

Further results in infinite dimensional topology are needed
to handle the nonregular points.

\section{The third order case}

In the second order periodic case, the set $C^\ast_{2,p}$
plays a very important role.
For the harder third order case, it is natural to start
by considering its counterpart,
the set $C^\ast_{3,p} \subset (H^3(\Ss^1))^2$
of pairs of potentials $(h_0, h_1)$
for which {\it all} solutions $v$ of 
\[   v'''(t) - h_1(t) v'(t) - h_0(t) v(t) = 0 \eqno{(\dagger)}\]
are periodic, i.e., satisfy
\[ v(0) = v(2\pi), \quad v'(0) = v'(2\pi), \quad v''(0) = v''(2\pi). \]
For $(h_0, h_1) \in C^\ast_{3,p}$,
let $v_0, v_1, v_2$ be the (fundamental) solutions
of $(\dagger)$ with initial conditions
\[ \begin{pmatrix} v_0 & v_1 & v_2 \\ v_0' & v_1' & v_2' \\
v_0'' & v_1'' & v_2'' \end{pmatrix}(0) = I \]
and normalize:
\[\gamma(t) = \frac{1}{|(v_0(t), v_1(t), v_2(t))|}\;(v_0(t), v_1(t), v_2(t)).\]
A straightforward computation yields
\[ \det(\gamma(0), \gamma'(0), \gamma''(0)) = 1, \qquad
\det(\gamma(t), \gamma'(t), \gamma''(t)) > 0 \textrm{ for all }t. \]
Conversely, a curve $\gamma: [0, 2\pi] \to \Ss^2$ is
{\it locally convex} if
$\det(\gamma(t), \gamma'(t), \gamma''(t)) > 0$ for all $t$.
Notice that this implies $\gamma'(t) \ne 0$ for all $t$.
Let $X_I$ be the set of locally convex curves $\gamma$
(with appropriate smoothness hypothesis) satisfying 
\[ \gamma(0) = \gamma(2\pi) = e_1, \quad \gamma'(0) = \gamma'(2\pi) = e_2,
\quad \gamma^{(j)}(0) = \gamma^{(j)}(2\pi) \]
and let $X_I^1 \subset X_I$ be the set of curves $\gamma$ with
$\det(\gamma(0), \gamma'(0), \gamma''(0)) = 1$.
Clearly, the inclusion $X_I^1 \subset X_I$ is a homotopy equivalence.
The map just described from $C^\ast_{3,p}$ to $X_I^1$
is a diffeomorphism; we proceed to construct its inverse.
Given $\gamma \in X_I^1$, set
\[ r(t) = \left(\det(\gamma(t), \gamma'(t), \gamma''(t))\right)^{-1/3} \]
and $V(t) = r(t) \gamma(t)$.
We then have $\det(V(t), V'(t), V''(t)) = 1$ for all $t$,
which implies that the vector $V'''(t)$
is a linear combination of $V(t)$ and $V'(t)$.
In other words, there exist unique real valued functions $h_0$ and $h_1$
with $V'''(t) = h_0(t) V(t) + h_1(t) V'(t)$, so that each coordinate
of $V$ is a periodic (fundamental) solution of $(\dagger)$
and the pair $(h_0,h_1)$ belongs to $C^\ast_{3,p}$.
Thus, $C^\ast_{3,p}$ is homotopy equivalent to $X_I$.

The space $X_I$ has been studied, among others, by
Little (\cite{Little}), B. Shapiro, M. Shapiro and Khesin
(\cite{Shapiro2}, \cite{ShapiroM}, \cite{SK}).
Little showed that $X_I$ has three connected components
which we shall call $X_{-1,c}$, $X_1$ and $X_{-1}$.
One of the authors (S., \cite{S}) established several
results concerning the homotopy and cohomology of these components.
It turns out that $X_{-1,c}$ is contractible and
$X_1$ and $X_{-1}$ are simply connected
but not homotopically equivalent to finite CW-complexes.
Also, $\pi_2(X_{-1})$ contains a copy of $\ZZ$
and $\pi_2(X_1)$ contains a copy of $\ZZ^2$.
As for the cohomology,
$H^n(X_{1}, \RR)$ and $H^n(X_{-1}, \RR)$
are nontrivial for all even $n$,
$\dim H^{4n-2}(X_1, \RR) \ge 2$ and
$\dim H^{4n}(X_{-1}, \RR) \ge 2$ for all positive $n$.

These results indicate that the topology of the critical set
of third order operators is far more complicated than
the lower order counterparts.

\bigskip\bigskip\bigbreak

{

\parindent=0pt
\parskip=0pt
\obeylines

Nicolau C. Saldanha and Carlos Tomei, Departamento de Matem\'atica, PUC-Rio
R. Marqu\^es de S. Vicente 225, Rio de Janeiro, RJ 22453-900, Brazil

\smallskip

nicolau@mat.puc-rio.br; http://www.mat.puc-rio.br/$\sim$nicolau/
tomei@mat.puc-rio.br

}

\end{document}